\documentclass[12pt,psamsfonts,openany]{article}
\usepackage{a4,amsmath,amssymb,amsopn,amsthm,graphics,latexsym}
\newtheorem{theorem}{Theorem}
\newtheorem{lemma}{Lemma}

\newtheorem{definition}{Definition}
\newtheorem{remark}{Remark}
\newtheorem{property}{Property}

\newcommand\me{\kern.3mm\mathbb{E}\kern.3mm}
\newcommand\dv{\kern.3mm{\rm{:}}\ }
\newcommand\tz{\kern.3mm{\rm{;}}\ }
\newcommand\sdots{.\kern.3mm.\kern.3mm.}

\newcommand {\beq}{\begin{equation}}
\newcommand {\eeq}{\end{equation}}
\newcommand {\barr}{\begin{eqnarray}}
\newcommand {\barn}{\begin{eqnarray*}}
\newcommand {\earr}{\end{eqnarray}}
\newcommand {\earn}{\end{eqnarray*}}
\newcommand{\bexam}{\begin{example}}
\newcommand{\eexam}{\end{example}}
\newcommand{\bpro}{\begin{property}}
\newcommand{\epro}{\end{property}}
\newcommand{\ble}{\begin{lemma}}
\newcommand{\ele}{\end{lemma}}
\newcommand{\bre}{\begin{remark}}
\newcommand{\ere}{\end{remark}}
\newcommand{\bdef}{\begin{definition}}
\newcommand{\edf}{\end{definition}}
\newcommand{\bth}{\begin{theorem}}
\newcommand{\enth}{\end{theorem}}
\newcommand{\bco}{\begin{corollary}}
\newcommand{\eco}{\end{corollary}}
\def\@{\nobreak\hskip0pt\hbox{\kern.5pt-}\hskip0pt\relax}

\begin{document}
\begin{center}
{\large\bf A note on exponential inequalities\break
for the distribution tails  of canonical Von Mises' statistics\\
of dependent observations}
\footnote{
Supported by the Russian Foundation for Basic Research,  grants
13--01--12415-OFI-M and 13--01--00511.}

\end{center}
\begin{center}
I. S. Borisov,  N. V. Volodko

\vspace{0.3 cm}

{\it Sobolev Institute of Mathematics, Novosibirsk State University,\\ Novosibirsk, 630090 Russia.\\
E-mail: sibam@math.nsc.ru, nvolodko@gmail.com}

\end{center}

%\maketitle

\begin{abstract}
\noindent The H{\"o}ffding-type exponential inequalities are obtained for the
distribution tails of {canonical}
%$U$- and
Von Mises' statistics of arbitrary order based on samples from a
stationary sequence of random variables satisfying~$\varphi$\@mixing condition.

\noindent{\it Keywords and phrases}: stationary sequence of random variables,
$\varphi$\@mixing, multiple orthogonal series, canonical $U$- and
$V$-statistics, exponential inequality, distribution tail, H{\"o}ffding's inequality.
\end{abstract}

\begin{center}
{\bf 1. Introduction. Preliminary results}
\end{center}

The paper deals with estimation of the distribution tails
of $U$- and $V$-statistics with canonical bounded kernels,
%based
%on samples of stationarily connected observations under~$\varphi$-mixing.
based on $n$ observations from a stationary $\varphi$-mixing process.
The exponential inequalities obtained are a natural generalization of
well-known H{\"o}ffding's inequality for the distribution tail of
a sum of independent identically distributed bounded random
variables. The approach of the present paper is  based on the kernel
representation of the statistics under consideration as a multiple
orthogonal series (for detail, see Borisov and Volodko, 2008, Korolyuk and
Borovskikh, 1994). It allows to
reduce the problem to more traditional estimates for the distribution tail of
a sum of weakly dependent random variables.

Introduce basic definitions and notions.
Let $X_1,X_2,\dots$~be a stationary sequence of random variables
taking values in an arbitrary measurable space~$\{\mathfrak X,\cal
A\}$ and having a common distribution~$F$. In addition to the stationary
sequence introduced above, we need an auxiliary
sequence~$\{X_i^{*}\}$ consisting of independent copies  of~$X_1$.

Given a natural $m$, denote by ~$L_2(\mathfrak{X}^m,F^m)$ the space of all measurable
functions~$f(t_1,\dots,t_m)$ defined on the corresponding
Cartesian power of the space~$\{\mathfrak X,\cal A\}$ with the
corresponding product measure and satisfying the condition $$
\mathbb{E}f^2(X_1^*, \ldots,X_m^*)<\infty. $$

\begin{definition}
{\rm A function $f(t_1,\dots,t_m)\in L_2(\mathfrak{X}^m,F^m)$ is called}
{\it canonical\/} if
%$$
\begin{equation}\label{deg}                       %(1)
\me f(t_1,...,t_{k-1},X_1,t_{k+1},...,t_m)=0
\end{equation}
%$$
%\vskip-7pt\noindent
{\rm for every $k=1,\ldots,m$ and all $t_j\in\mathfrak{X}$ under obvious agreement of the notation in the cases $k=1$ and $k=m$.}
%conditional expectation given the random variables $\{X^*_i; \,i\le m, i\neq k\}$.
\end{definition}
Define a Von Mises' statistic (or $V$-{\it statistic}) by the
formula
%\vskip-0.001pt\noindent
%$$
\begin{equation}\label{Mis}                         %(2)
V_n\equiv V_n(f):=n^{-m/2} \sum_{1\leq j_1,\dots,j_m\leq
n}f(X_{j_1},\dots,X_{j_m}).
\end{equation}
In the sequel, we consider only the statistics where the function
$f(t_1, \ldots, t_m)$ (the so-called \textit{kernel} of the
statistic) is canonical. In this case, the corresponding $V$-statistic is also called \textit{canonical}. For
independent~$\{X_i\}$, these statistics are studied during last
sixty years. The corresponding reference and examples of such
statistics can be found in Korolyuk and Borovskikh, 1994. Notice also that any statistic
having the structure of the Euclidean norm squared of the $n$th partial normed sum of random variables
(independent or not) taking values in a Hilbert space, can be represented in the form
(\ref{Mis}). For example, the classical $\omega^2$-statistics, $\chi^2$-statistics, and some others have the structure mentioned.

In addition to $V$-statistics, the
so-called $U$-statistics are studied as well:
\begin{equation}\label{Ustat}
U_n\equiv U_n(f):=\left (\frac{(n-m)!}{n!}\right )^{1/2} \sum \limits_{1 \leq i_1 \neq
\dots  \neq i_m \leq n} f(X_{i_1}, \ldots , X_{i_m}).
\end{equation}
In this case, the value $(n-m)!/n!$ is equivalent to $n^{-m}$ as $n\to\infty$.
Notice also that any $U$-statistic is represented as a finite
linear combination of canonical $U$-statistics of the orders from~1
to~$m$. This representation is called  {\it H{\"o}ffding's decomposition\/}
(see~Korolyuk and Borovskikh, 1994). Every $V$-statistic ($U$-statistic) can be represented as a linear combination of $U$-statistics ($V$-statistic, respectively) of the orders from $1$ to $m$.
So, to estimate the distribution tails of canonical statistics from one of these two classes, we may estimate those for the second one.

For independent observations~$\{X_i\}$, we give below a brief review
of the results directly connected with the topic of the present paper.
%In
%this connection, we would like to mention the results in \cite[Theorem~1]{B},
%\cite[Proposition~2.2]{Gine}, \cite[Theorem~7, Corollary~3]{Adam},
%and \cite[Theorem~3.3]{GLZ}.
One of the first papers dealing with exponential inequalities for the distribution tails of
$U$-statistics was the paper by
H{\"o}ffding, 1963, although he considered noncanonical $U$-statistics only. %In this case, the value $(n-m)!/n!$
%equivalent to $n^{-m}$ as $n\to\infty$, is used as the normalizing
%factor instead of $n^{-m/2}$.
In particular, in the case $m=1$, the following
statement is contained in H{\"o}ffding, 1963 as a particular case:
\begin{equation}\label{H}
\mathbb{P}( U_n-\mathbb{E} U_n\geq t)\leq e^{-2t^2/(b-a)^2},
\end{equation}
%\vskip-7pt\noindent
where
%\vskip-1pt\noindent
%$$ \tilde U_n:=(n-m)!/n!\sum \limits_{1 \leq i_1 \neq
%\dots  \neq i_m \leq n} f(X_{i_1}, \dots , X_{i_m}), $$
$a\leq f(t_1, \ldots,t_m)\leq b$.
%In the case $m=1$,
Inequality (\ref{H}) is usually called {\it H{\"o}ffding's
inequality\/} for sums of independent identically distributed
bounded random variables. Notice that, in this case, the
centered sums mentioned may be considered as  the simplest example of canonical
%or  noncanonical
$V$-statistics.

In Borisov, 1990, 1991, an improvement of (\ref{H}) was obtained for the
case when there exists a splitting majorant of the canonical
kernel under consideration:
\begin{equation}\label{M}
\big|f(t_1,\dots ,t_m)\big|\le \prod_{i\le m}g(t_i),
\end{equation}
and the function $g(t)$ satisfies Bernstein's condition $$ {\mathbb
E}g(X_1)^k\le \sigma^2L^{k-2}k!/2 $$
for all $k\ge 2$, where $\sigma$ and $L$ are some positive constants. In this
case, the following analogue of Bernstein's inequality holds:
\begin{equation}\label{Bor}
\mathbb{P}\big(|V_n|\geq t\big)\leq
2\exp\left(-\frac{c_1t^{2/m}}{\sigma^2+Lt^{1/m}n^{-1/2}}\right),
\end{equation}
where the constant $c_1$ depends  on~$m$ only. Moreover,
%as noted in~\cite{B0} and \cite{B},
inequality~(\ref{Bor}) cannot be
improved in a sense.

It is clear that if $\sup_{t_i}\big|f(t_1, \dots,
t_m)\big|=B<\infty$ then one can put $g(\cdot)=\sigma=L=B^{1/m}$   in (\ref{M}) and (\ref{Bor}). Then it
suffices to consider only the deviation range $|t|\le
Bn^{m/2}$ in~(\ref{Bor}) (otherwise, the left-hand side of~(\ref{Bor})
vanishes). Therefore, for all $t\ge 0$, inequality (\ref{Bor}) yields the
upper bound

\begin{equation}\label{Arcones}
\mathbb{P}\big(|V_n|\geq t\big)\leq
2\exp\left(-\frac{c_1}{2}(t/B)^{2/m}\right)
\end{equation}
which is an analogue of H{\"o}ffding's inequality~(\ref{H}) for the $m$th
power of the normed sum of independent identically distributed centered and bounded random variables,
i.~e., for an elementary example of canonical $V$-statistics of order $m$.

In Arcones and Gine, 1993, some inequality close to~(\ref{Bor})  was proved
without condition~(\ref{M}), and relation~(\ref{Arcones}) is given
as a consequence. In Gine et al., 2000, some refinement of (\ref{Arcones})
is obtained for $m=2$, and in Adamczak, 2006, the later result was extended
to canonical $U$-statistics of an~arbitrary
order. The most accurate estimates in the later case was obtained in Major, 2006 
and Major, 2007.

An extension of inequality (\ref{Arcones}) to the case of
stationarily connected observations under $\varphi$-mixing condition is adduced in
Borisov and Volodko, 2009a, 2009b.
The goal of the present paper is to weaken the restrictions on the coefficient $\varphi(\cdot)$
which are contained in Borisov and Volodko, 2009a, 2009b.
At the same time,  we correct  the corresponding proof in
 the later papers.
\pagebreak

\begin{center}
{\bf 2. Main results for weakly dependent observations}
\end{center}

In the sequel, we assume that $\mathfrak{X}$ is a separable metric space. Then the Hilbert
space~$L_2(\mathfrak{X},F)$ has a countable orthonormal
basis~$\{e_i(t)\}$. Put~$e_0(t)\equiv 1$.
Using the Gram--Schmidt orthogonalization, one can construct an~orthonormal basis
in~$L_2(\mathfrak{X},F)$ containing the constant function
$e_0(t)$. Then $\mathbb{E}e_i(X_1)=0$ for every $i\geq 1$ due to
orthogonality of all the other basis elements to the
function~$e_0(t)$. The normalizing condition means that~$\me
e^2_i(X_1)=1$ for all~$i\geq 1$.
We assume that the basis consists of uniformly bounded functions:
%\vskip-7pt\noindent
\begin{equation}\label{boundbase}
\sup_{i, t}|e_i(t)|\leq C.
\end{equation}

It is well-known that the collection of the functions $$
\big\{e_{i_1}(t_1)e_{i_2}(t_2)\cdots e_{i_m}(t_m); \quad i_2,
\dots, i_m=0, 1, \dots\big\} $$ is an orthonormal basis of the
Hilbert space $L_2(\mathfrak{X}^m,F^m)$. The kernel $f(t_1,\dots ,t_m)$
can be decomposed by the basis $\big\{e_{i_1}(t_1)\cdots
e_{i_m}(t_m)\big\}$ and represented as the series
%\vskip-3pt\noindent
\begin{equation}\label{kernel}
f(t_1,\dots ,t_m)=\sum_{i_1,\dots ,i_m=1}^{\infty} f_{i_1,\dots,
i_m}e_{i_1}(t_1)\cdots e_{i_m}(t_m)
\end{equation}
which converges in the norm of
$L_2(\mathfrak{X}^m,F^m)$. The basis element $e_0(t)$ is absent in
representation~(\ref{kernel}) because the kernel is canonical
(for detail, see Borisov and Volodko, 2008). Moreover, if the
coefficients~$\{f_{i_1,\dots,i_m}\}$ are absolutely summable then,
due to the B. Levi theorem and the simple estimate ${\bf
E}\big\vert e_{i_1}(X^*_1)\cdots\,e_{i_m}(X^*_m)\big\vert\le 1$,
the series in~(\ref{kernel}) converges almost surely with respect
to the distribution $F^m$ of the vector $(X^*_1,\dots,X^*_m)$. It
is worth noting that
%, even in this case,
we cannot extend the last claim
to the case when the  vector $(X_1,\dots ,X_m)$ has dependent
coordinates. If we substitute dependent random variables $X_1,...,X_m$ for the nonrandom argument $t_1,...,t_m$ in (\ref{kernel}) then
equality (\ref{kernel}) may be in general false with a nonzero probability (see the
corresponding example in Borisov and Volodko, 2008). It is explained by the fact 
that
the exclusive sets for the distributions of the random vectors $(X_1,...,X_m)$
and $(X^*_1,...,X^*_m)$ may essentially differ.
There exist different ways to avoid this difficulty, for example,
one may require absolute continuity of the distribution of the random vector
$(X_1,...,X_m)$ with respect to that of
$(X^*_1,...,X^*_m)$ or continuity both of the kernel and basis elements
(see Borisov and Volodko, 2008).
If the kernel $f(\cdot)$ is defined by
%one considers the right-hand side of (\ref{kernel}) as $f(\cdot)$, i.~e.
%requires the fulfilment of
equality (\ref{kernel}) everywhere, i.~e., for all values of the vector
argument, then this equality will be naturally valid  after replacement of the vector argument by
arbitrarily connected observations $X_1,...,X_m$ in (\ref{kernel}) and it will be no
problems in this case.

Now, assume that after the above-mentioned replacement of the nonrandom argument in  (\ref{kernel}) we preserve equality (\ref{kernel}) almost surely and substitute the resulting relation into (\ref{Mis}). Then the following key representation is valid:
%the definition of the Von Mises statistic
\begin{align*}
V_n &=n^{-m/2}\sum_{1\leq j_1, \dots, j_m\leq n}f(X_{j_1}, \dots,
X_{j_m})\\
    &=n^{-m/2}\sum_{1\leq j_1, \dots, j_m\leq n}
      \sum_{i_1,\dots ,i_m=1}^{\infty}f_{i_1,\ldots, i_m}e_{i_1}(X_{j_1})
       \cdots e_{i_m}(X_{j_m})\\
    &=\sum_{i_1,\dots ,i_m=1}^{\infty}f_{i_1,\ldots, i_m}n^{-1/2}
      \sum_{j=1}^{n}e_{i_1}(X_{j})\cdots
      n^{-1/2}\sum_{j=1}^{n}e_{i_m}(X_{j})\\
    &=\sum_{i_1,\dots ,i_m=1}^{\infty}f_{i_1,\ldots, i_m}
      S_n(i_1)\cdots S_n(i_m),
\end{align*}
where    $S_n(i_k):=n^{-1/2}\sum_{j=1}^{n}e_{i_k}(X_{j})$,\,\,$k=1,...,m$.

%$$
\goodbreak

In the present paper, we consider only stationary
sequences~$\{X_j\}$ satisfying  $\varphi$-mixing condition.
Recall the definition of this type of dependence. For~$j\le k$,
denote by~$\mathfrak{M}^k_{j}$ the $\sigma$-field of all events
generated by the random variables $X_j,\dots,X_k$.

\begin{definition}
%\bdef
{\rm A  sequence $X_1,\ X_2,\dots$ satisfies $\varphi$-{\it
mixing\/} (or {\it uniformly strong mixing}) if $$
\varphi(i):=\sup_{k\ge 1}\ \sup_{A\in\mathfrak{M}^k_{1},\,B\in
\mathfrak{M}^{\infty}_{k+i},\,\mathbb{P}(A)>0}
\frac{|\mathbb{P}(AB)-\mathbb{P}(A)\mathbb{P}(B)|}
{\mathbb{P}(A)}\rightarrow 0\quad \mbox{as}\,\,\,i\rightarrow\infty. $$}
%\edf
\end{definition}

\bre
{\rm If $\underline{}\{X_j\}$ satisfies $\varphi$-mixing with
coefficient~$\varphi(\cdot)$ then, for every measurable function~$f$, the sequence $f(X_1),f(X_2),\dots $
 also satisfies $\varphi$-mixing
with some coefficient which does not exceed~$\varphi(\cdot)$.
Notice also that we consider $\varphi(\cdot)$ as a function defined on the set of all nonnegative integers.
%${\mathbb Z}_+:=\{i:\,i=0, 1,\ldots\}$.
It is clear that $\varphi(0)=1$ for an arbitrary stationary sequence.}
\ere

The main result of the present paper  is as follows.

%\bth
{\bf Theorem}. {\it Let a canonical kernel~$f(t_1,\dots,t_m)$ and the basis function $\{e_k(t)\}$
be continuous $($in the product topology$)$ everywhere on~$\mathfrak X^m$ and let condition $(\ref{boundbase})$ be fulfilled.
Moreover, if
$$\sum_{i_{1},\dots ,i_{m}=1}^{\infty}|f_{i_{1},\ldots, i_{m}}|<\infty$$
and
\begin{equation}\label{phispeed}
 \Phi:=\sum_{k=0}^{\infty} \varphi(k)< \infty
\end{equation}
then the following inequality
holds$:$
\begin{equation*}
%\label{first}
{\mathbb P}\big(|V_n|>x\big)\leq
\exp\big\{-(16\Phi e)^{-1}(x/B(f))^{2/m}\big\},
\end{equation*}
where
%$C_1>0$ depends only on~$\Phi$,
$ B(f):=C^m\sum_{i_{1},\dots ,i_{m}=1}^{\infty}|f_{i_{1},
\ldots, i_{m}}|$ and the constant $C$ is defined in $(\ref{boundbase})$.}
%\enth

\bre
{\rm If the kernel $f(\cdot)$ can be represented as (\ref{kernel}) for {\it all values}
of the vector argument, then the requirement of continuity of the kernel and
the basis functions in the theorem conditions is unnecessary. Moreover, the requirement
of continuity can be omitted if the distribution of the random vector
$(X_1,...,X_m)$ is absolute continuous with respect to that of $(X^*_1,...,X^*_m)$ (see Borisov and Volodko, 2008).}
\ere

\bre
{\rm The proof of the theorem above is much shorter than the corresponding
proofs from Borisov and Volodko, 2009a, 2009b, where an exponential decreasing of
the coefficient $\varphi$ was required. This restriction is much stronger than
(\ref{phispeed}). Moreover, the proof from Borisov and Volodko, 2009a contains essential
inaccuracy. The crucial point of the present proof is a moment inequality for
sums of stationarily connected random variables under $\varphi$-mixing, obtained in
Dedecker and Prieur, 2005.}
\ere

%\begin{center}
{\it Proof of Theorem}.
%\end{center}
%\begin{proof}
Without loss of generality, we assume that the separable metric
space~$\mathfrak X$ coincides with the support of the
distribution~$F$. The last means that $\mathfrak X$ does not
contain the open balls with $F$-measure zero. Since all the
basis elements $e_k(t)$ in~(\ref{kernel}) are continuous and
uniformly bounded in~$t$ and $k$, due to  Lebesgue's dominated
convergence theorem, the series in~(\ref{kernel}) is continuous if
the coefficients $f_{i_1,\ldots, i_m}$ are absolutely summable. It
is not difficult to see that, in this case,  the equality
in~(\ref{kernel}) turns into the identity on the all variables
$t_1,\dots ,t_m$ because  equality of two continuous functions
on an everywhere dense set implies  their coincidence everywhere.
So, in this case,  one can substitute
{\it arbitrarily dependent\/} observations for the variables
$t_1,\dots ,t_m$ in identity~(\ref{kernel}). Therefore, for {\it all elementary events\/}, the
above-mentioned representation holds:
\begin{equation}\label{multser}
V_n=\sum_{i_1,\ldots,i_m=1}^{\infty}f_{i_1,\ldots, i_m}S_n(i_1)\cdots  S_n(i_m),
\end{equation}
where, as above,    $S_n(i_k):=n^{-1/2}\sum_{j=1}^{n}e_{i_k}(X_{j})$,\,\,$k=1,...,m$.

%\pagebreak
Consider an arbitrary even moment
of the above-introduced $V$-statistic. First, from (\ref{multser}) we have
%$$
\begin{eqnarray}\label{NMis}
\mathbb{E}V_n^{2N}=\sum_{i_1,\dots
,i_{2mN}=1}^{\infty} f_{i_1,\ldots, i_m}\cdots
f_{i_{2mN-m+1},\ldots, i_{2mN}}\mathbb{E}S_n(i_1) \cdots
S_n(i_{2mN})\nonumber \\
\leq \sum_{i_1,...,i_{2mN}=1}^{\infty}
|f_{i_1...i_m}|...|f_{i_{2mN-m+1}...i_{2mN}}|(\mathbb{E}S_n^{2mN}(i_1))^{1/2mN}...
(\mathbb{E}S_n^{2mN}(i_{2mN}))^{1/2mN}.\nonumber \\
%\mbox{where}\quad N=0, 1, 2,\ldots.\qquad\qquad\qquad\qquad\qquad\qquad\qquad\qquad\qquad\qquad\qquad\qquad
\end{eqnarray}
Next we need the following auxiliary statement from Proposition 5 of Dedecker and Prieur, 2005,
adapted to our conditions.
%\bth

{\bf Proposition}.
{\it Let $Y_1,Y_2,...$ be a stationary sequence
of random variables taking values in $\mathbb{R}$ and satisfying  $\varphi$-mixing condition , moreover, $|Y_1|\leq C$ almost surely.
% $h(t):\mathbb{R}\rightarrow\mathbb{R}$ --- дг­ЄжЁп б ®Ја ­ЁзҐ­­®© ў аЁ жЁҐ©, в.Ґ. Є®­Ґз­  б«Ґ¤гой п ўҐ«ЁзЁ­ :
% $$\|dh\|=\sup_P\sum_{k=0}^m|h(x_k+1)-h(x_k)|,$$
% Ј¤Ґ $P=(x_0,...,x_{m+1})$ --- Їа®Ё§ў®«м­®Ґ а §ЎЁҐ­ЁҐ ®Ў« бвЁ ®ЇаҐ¤Ґ«Ґ­Ёп $h$. Ља®¬Ґ в®Ј®,
% Є®­Ґз­® $\mathbb{E}h(Y_1)$.
Then, for every $p\geq 2$, the following inequality is valid:}
\begin{equation}\label{exp}
\mathbb{E}|\sum_{i=1}^{n}Y_i-n\mathbb{E}Y_1|^p\leq\big(8C^2p\sum_{k=0}^{n-1}(n-k)
\varphi(k)\big)^{p/2}.
\end{equation}

Now, for every fixed $i$, consider the sequence
$e_i(X_1), e_i(X_2),...$  as the sequence of random variables $\{Y_j\}$ in (\ref{exp}). Then we can obtain the following estimate for the even moments of the above-introduced normed sums: $$\mathbb{E}S_n^{2mN}(i)\leq \big(16\Phi C^2mN\big)^{mN},\qquad N=0,1,\ldots ,$$
%\sum_{k=0}^{\infty}\varphi(k)\big)^{mN},$$
and from this and (\ref{NMis}) we deduce the corresponding estimate for the even  moments of the $V$-statistic:
$$\mathbb{E}V_n^{2N}\leq\Big(C^m\sum_{i_1,...,i_m=1}^{\infty}|f_{i_1...i_m}|\Big)^{2N}
(16\Phi mN)^{mN}.$$
%where $\tilde{c}:=$

Now we can estimate the distribution tail of the $V$-statistic using Chebyshev's inequality as follows:
$${\mathbb P}(|V_n|>x)\leq{x^{-2N}}{\mathbb{E}V_n^{2N}}\leq x^{-2N}{\Big(C^m\sum_{i_1,...,i_m=1}^{
\infty}|f_{i_1...i_m}|\Big)^{2N}(16\Phi mN)^{mN}}.$$
Put $N=[\varepsilon x^{2/m}]$, where $\varepsilon>0$ is arbitrary and $[a]$ is the integer part of a positive number $a$.
Then we have
$${\mathbb P}(|V_n|>x)\leq {x^{-2N}{\Big(C^m\sum_{i_1,...,i_m=1}^{
\infty}|f_{i_1...i_m}|\Big)^{2N}(16\Phi m\varepsilon)^{mN}x^{2N}}}$$
$$\le\exp\{\varepsilon mx^{2/m}\log(\tilde{c}m\varepsilon)\},$$
where $$\tilde{c}:=
%16\Phi\Big(C^m\sum_{i_1,...,i_m=1}^{\infty}|f_{i_1...i_m}|\Big)^{2/m}=
16\Phi (B(f))^{2/m}.$$
The multiplier $\varepsilon m\log(\tilde{c}m\varepsilon)$
attains its minimal value at the point $\varepsilon_0:=(\tilde{c}me)^{-1}$ and this
minimal value equals $-(\tilde{c} e)^{-1}$. Thus,
$${\mathbb P}(|V_n|>x)\leq \exp\{\varepsilon_0 mx^{2/m}\log(\tilde{c}m\varepsilon_0)\}=\exp\{-(16\Phi e)^{-1}(x/B(f))^{2/m}\}$$ which
was to be proved.
%\end{proof}

\end{document}